\documentclass[12pt,a4paper,draft]{article}
\usepackage{amssymb,amsmath}
\begin{document}

\begin{flushright}
UDC 512.542
\end{flushright}

\begin{center}
{\bfseries\scshape ON THE CONJUGACY PROBLEM FOR CARTER SUBGROUPS}
\bigskip

E.P. Vdovin\footnote{The authors were supported by the~Russian
Foundation for Basic Research (Grant 05--01--00797),
the Program ``Universities of Russia'' (Grant UR.04.01.202), and by Presidium
SB RAS Grant 86-197.}
\end{center}

\begin{center}
{\bfseries Abstract}
\end{center}
{\small It is proven in the paper, that Carter subgroups of a finite group are 
conjugate if Carter subgroups in the group of induced automorphisms for every
non-Abelian composition factor are conjugate.}
\vspace{0.5\baselineskip}

{\bfseries Key words.} Carter subgroup, almost simple group, group of induced
automorphisms.

\begin{center}
{\bfseries\Large \S~1. Introduction}
\end{center}

Recall that a nilpotent self-normalizing subgroup is called a {\em Carter
subgroup}. In the paper we consider the following
\vspace{0.5\baselineskip}

\noindent{\bfseries\itshape Problem.} Are any two Carter subgroups of a finite
group conjugate?
\vspace{0.5\baselineskip}

In [1] it is proven that the minimal counter example to this
problem should be almost simple.  We intend to improve results obtained in
[1] (see the theorem below). Actually, we shall use ideas of [1]
in order to
prove stronger theorem.

Our notations is standard.    For a finite
group $G$ we denote by $\mathrm{Aut}(G)$ the
group of automorphisms of $G$. If $Z(G)$ is trivial, then $G$ is
isomorphic to the group of its inner automorhisms and we may suppose
that~$G\leq \mathrm{Aut}(G)$.  A finite group $G$ is said to be {\em almost
simple} if there is a simple group $S$ with $S\leq G\leq \mathrm{Aut}( S)$,
i.~e., $F^\ast(G)$ is a simple group. We denote by $F(G)$ the Fitting subgroup
of $G$ and by $F^\ast(G)$ the generalized Fitting subgroup of~$G$.

If $G$ is a group, $A,B,H$ are subgroups of $G$ and $B$ is normal in $A$
($B\unlhd A$), then $N_H(A/B)=N_H(A)\cap N_H(B)$. If $x\in N_H(A/B)$, then $x$
induces an automorphism $Ba\mapsto B x^{-1}ax$ of $A/B$. Thus, there is a
homomorphism of $N_H(A/B)$ into $\mathrm{Aut}(A/B)$. The image of this
homomorphism is
denoted by $\mathrm{Aut}_H(A/B)$ while its kernel is denoted by $C_H(A/B)$. In
particular, if $S$ is a composition factor of $G$, then for any $H\leq G$ the
group $\mathrm{Aut}_H(S)$ is defined. 
\vspace{0.5\baselineskip}

\noindent{\bfseries Definition.} 
A finite group $G$ is said to satisfy condition
($*$) if, for every its non-abelian composition factor $S$ and
for every its nilpotent subgroup $N$, Carter subgroups of $\langle
\mathrm{Aut}_N(S),S\rangle$ are conjugate.
\vspace{0.5\baselineskip}
 
Clearly,  if a finite group $G$ satisfies ($*$), then for every normal
subgroup $H$  and soluble subgroup $N$ of $G$, groups $G/H$ and $NH$
satisfy ($*$). Our goal here is to prove the following theorem.
 \vspace{0.5\baselineskip}
 
 \noindent{\bfseries Theorem.}
{\itshape If a finite group $G$ satisfies {\em ($*$)}, then Carter subgroups of
$G$ are conjugate.}
 \vspace{0.5\baselineskip}

Note that a finite group may not contain Carter subgroups. In this case we also
say that its Carter subgroups are conjugate. In sections 2 and 3 we are
assuming that $X$ is a counter example to the theorem of minimal order, i.~e.,
that $X$ is a finite group satisfying condition ($*$), $X$ contains
nonconjugate Carter subgroups, but Carter subgroups in every group $M$  of order
less, that $\vert X\vert$, satisfying condition ($*$), are conjugate.

\begin{center}
{\bfseries\Large \S~2. Preliminary results}
\end{center}

Recall that $X$ is a counter example to the theorem of minimal order.
\vspace{0.5\baselineskip}

\noindent{\bfseries Lemma 1.} {\itshape
Let $G$ be a finite group satisfying {\em ($*$)} and $\vert G\vert\leqslant
\vert X \vert$. Let  $H$ be a Carter subgroup of $G$. If $N$ is a normal
subgroup of $G$, then $HN/N$ is a Carter subgroup of~${G/N}$.}
\vspace{0.5\baselineskip}

\noindent{\itshape Proof.}
Since $HN/N$ is nilpotent, we have just to prove that it is self-normalizing in
$G/N$. Clearly, this is true if $G=HN$. So, assume $M=HN<G$.  By the minimality
of $X$, $M^x=M$, $x\in G$, implies $H^x=H^m$ for some $m\in M$. It follows
$xm^{-1}\in N_G(H)=H$ and $x\in M$. This proves that $HN/N$ is nilpotent and
self-normalizing in~$G/N$. \hfill{$\Box$}
\vspace{0.5\baselineskip}

\noindent{\bfseries Lemma 2.} {\itshape
Let $B$ be a minimal normal subgroup of $X$ and $H,K$ be non-conjugate Carter
subgroups of $X$.
\begin{itemize}
\item[{\em (i)}] $B$ is non-soluble;
\item[{\em (ii)}] $X=BH=BK$;
\item[{\em (iii)}] $B$ is the unique minimal normal subgroup of~$X$.
\end{itemize} }
\vspace{0.5\baselineskip}

\noindent{\itshape Proof.}
(i) We give a proof by contradiction. Assume that $B$ is soluble and let
$\pi:X\rightarrow X/B$ be the canonical homomorphism. Then $H^\pi$ and $K^\pi$
are Carter subgroups of $X/B$, by Lemma 1. By the
minimality of $X$, there  exists $\bar{x}=Bx$ such that
$(K^\pi)^{\bar{x}}=H^\pi$. It follows $K^x\leq BH$. Since $BH$ is soluble,
$K^x$ is conjugate to $H$ in $BH$, hence $K$ is conjugate to $H$ in $X$, a
contradiction.

(ii) Assume that $BH<X$. By Lemma 1 and the minimality of
$X$, $BH/B$ and $BK/B$ are conjugate in $X/B$: so there exists $x\in X$ such
that $K^x\leq BH$. It follows that $K^x$ is conjugate to $H$ in $BH$, hence $K$
is conjugate to $H$ in $X$, a contradiction.

(iii) Suppose that $M$ is a minimal normal subgroup of $X$ different from $B$.
By (i), $M$ is non-soluble. On the other hand, $MB/B\simeq M$ is a subgroup of
the nilpotent group $X/B\simeq H/H\cap B$, a contradiction. \hfill{$\Box$}
\vspace{0.5\baselineskip}

The following lemma is useful in many applications, so we prove it here, though
we need only a part of its proof in our later arguments.
\vspace{0.5\baselineskip}

\noindent{\bfseries Lemma 3.} {\itshape
Let $G$ be a finite group. Let $H$ be a Carter
subgroup of $G$. Assume that there exists a normal subgroup
$B=T_1\times\ldots\times T_k$ of $G$ such that $T_1\simeq\ldots\simeq T_k\simeq
T$, $Z(T_i)=\{1\}$ for all $i$, and $G=H(T_1\times\ldots\times T_k)$. Then  
$\mathrm{Aut}_H(T_i)$ is a Carter subgroup of~$\langle \mathrm{Aut}_H(T_i),
T_i\rangle$. }
\vspace{0.5\baselineskip}

\noindent{\itshape Proof.}
Assume that our statement is false and $G$ is a counterexample with $k$
minimal, then $k>1$. Clearly, $G$ acts
transitively, by conjugation, on the set $\Omega:=\{T_1,\ldots, T_k\}$. We may
assume that the $T_j$'s are indexed so that $G$ acts
primitively on the set $\{\Delta_1,\ldots,\Delta_p\}$, $p>1$, where for
each~$i$: $$\Delta_i:=\{T_{1+(i-1)l},\ldots,T_{il}\},\ \ \ \ \ k=pl.$$ 
Denote by $\varphi:G\rightarrow \mathrm{Sym}_p$  the induced permutation
representation. Clearly, $B\leq ker\ \varphi$, so that
$G^\varphi=(BH)^\varphi=H^\varphi$ is a primitive nilpotent subgroup of
$\mathrm{Sym}_p$. Hence $p$ is prime and $G^\varphi$ is a cyclic group of order
$p$. In
particular, $Y:=ker\ \varphi$ coincides with the stabilizer of any $\Delta_i$,
so that $\varphi$ is permutationally equivalent to the representation of $G$ on
the right cosets of $Y$. For each $i=1,\ldots,p$, let
$S_i=T_{1+(i-1)l}\times\ldots\times T_{il}$: then $Y=N_G(S_i)$ and
$B=S_1\times\ldots\times S_p$. Consider $\xi:Y\rightarrow \mathrm{Aut}_Y(S_1)$,
let
$A=Y^\xi$, $S=S_1^\xi$. Clearly $S$ is a normal subgroup of $A$; moreover, $S$
is isomorphic to $S_1$, since $S_1$ has trivial center. On the other hand, for
each $i\not=1$, $S_i\leq ker\ \xi$, since $S_i$ centralizes~$S_1$. 

Denote by $A\wr C_p$ the wreath product of $A$ and a cyclic group $C_p$ and let
$\{x_1=e,\ldots,x_p\}$ be a right transversal of $Y$. Then the map
$\eta:G\rightarrow A\wr C_p$ such that, for each $x\in G$:
$$x\mapsto \left(\left(x_1xx^{-1}_{1^{x^\varphi}}
\right)^\xi,\ldots,\left(x_pxx^{-1}_{p^{x^\varphi}}
\right)^\xi\right)x^\varphi$$ is a homomorphism. Clearly $Y^\eta$ is a
subdirect product of the base subgroup $A^p$ and $$S_1^\eta=
\{(s,1,\ldots,1)\vert s\in S\}, B^\eta=\{(s_1,\ldots,s_p)\vert s_i\in S\}\leq
Y^\eta.$$ Moreover, $ker\ \eta=C_G(B)=\{e\}$, so we may identify $G$ with
$G^\eta$. We choose $h\in H\setminus Y$. Then $G=\langle Y,h\rangle$, $h^p\in
Y$, $H=(Y\cap H)\langle h\rangle$ and we may assume
$$h=(a_1,a_2,\ldots,a_p)\pi, a_i\in A, \pi=(1,2,\ldots,p)\in C_p.$$ For each
$i$, $1\leqslant i\leqslant p$, let $\psi_i:A^p\rightarrow A$ be the canonical
projection
and let $H_i:=(H\cap Y)^{\psi_i}$. Clearly, $Y^{\psi_i}=A$. Moreover, for each
$i\geqslant 2$, $H_i=H_1^{h^{i-1}}=H_1^{a_1\ldots a_{i-1}}$ since $h$
normalizes~${Y\cap H}$. Let $N:=(H_1\times\ldots\times H_p)\cap Y.$ $N$ is
normalized by $H$, since $H=(N\cap H)\langle h\rangle$ and $H_i^h=H_{i+1\pmod
p}.$ We claim that  $H_1$ is a Carter subgroup of $A$. Assume $n_1\in
N_A(H_1)\setminus H_1$. From $Y=(Y\cap H)B$, it follows $n_1=h_1s$, $h_1\in
H_1$, $s\in N_S(H_1)\setminus H_1$. Let $b:=(s,s^{a_1},\ldots,s^{a_1\ldots
a_{p-1}})\in B.$ Then $b$ normalizes $N$, for:
$$H_i^b=H_i^{s^{a_1\ldots
a_{i-1}}}=H_1^{a_1\ldots a_{i-1}s^{a_1\ldots
a_{i-1}}}=H_1^{sa_1\ldots a_{i-1}}=H_1^{a_1\ldots a_{i-1}} =H_i.$$ Now
$[b,h^{-1}]:=b^{-1}hbh^{-1}\in Y$ is such that: $$[b,h^{-1}]^{\psi_i}=1\text{
if }i\not=p, [b,h^{-1}]^{\psi_p}=[s,(a_1\cdot\ldots\cdot
a_p)^{-1}]^{a_1\cdot\ldots\cdot
a_{p-1}},$$ where $a_1\cdot\ldots\cdot a_p=(h^p)^{\psi_1}\in H_1$. Since $s\in
N_S(H_1)$,
it follows $$[s,(a_1\cdot\ldots\cdot a_p)^{-1}]\in H_1,\
[s,(a_1\cdot\ldots\cdot a_p)^{-1}]^{a_1\cdot\ldots\cdot a_{p-1}}\in H_p.$$ So
$[b,h^{-1}]\in N$
and $b\in  N_G(N\langle h\rangle)$. But $H\leq N\langle h\rangle$, implies
$N_G(N\langle h\rangle)=N\langle h\rangle$. Indeed, if $g\in N_G(N\langle
h\rangle)$, then $H^g$ is a Carter subgroup of $N\langle h\rangle$. But
$N\langle h\rangle$ is soluble, hence there exists $y\in N\langle h\rangle$
with $H^g=H^y$. Now $H$ is a Carter subgroup of $G$, thus $gy^{-1}\in H$
and $g\in N\langle h\rangle$. Therefore $b\in N, s\in H_1$, i.~e., $n_1\in H_1$,
a contradiction.

Now $A=H_1(T_1\times\ldots\times T_l)$ and $l<k$. By induction we have that
$\mathrm{Aut}_{H_1}(T_1)$ is a Carter subgroup of
$\langle\mathrm{Aut}_{H_1}(T_1),T_1\rangle$.
In view of our construction, $\mathrm{Aut}_H(T_1)=\mathrm{Aut}_{H_1}(T_1)$ and
the lemma follows. \hfill{$\Box$}

\begin{center}
{\bfseries\Large \S~3. Proof of the Theorem}
\end{center}

Write $B=T_1\times\ldots\times T_k$, $T_i\simeq T$, a non-abelian simple group.
What remains to prove is $k=1$. In the notatinons of the proof of Lemma 3, we
have shown that $H_1$ is a Carter subgroup of
$A$. Clearly, since each $H_i$ is conjugate to $H_1$ in $A$, $N_A(H_i)=H_i$,
$i=1,\ldots,p$. It follows easily that $N$ is a Carter subgroup of $Y$. For,
let $y:=(y_1,\ldots,y_p)\in N_Y(N)$: from $N^{\psi_i}=H_i$, we have $y_i\in
N_A(H_i)=H_i$, for each $i$, hence~${y\in N}$.

We have seen that, to each Carter subgroup $H$ of $X$ we can associate a Carter
subgroup $N=N_H$ of $Y$, such that $H$ normalizes $N_H$. Clearly,
$N_H\not=\{e\}$, otherwise $X$ would have order $p$. So let $K$ be a Carter
subgroup of $X$, not conjugate to $H$, and let $N_K$ be the Carter subgroup of
$Y$ corresponding to $K$. By the minimality of $X$, for each $x\in X$ we have
$\langle H^x,K\rangle=X$. On the other hand, by the inductive hypothesis, there
exists $x\in Y$ such that $N_K=(N_H)^x$. Hence $N_K$ is normal in $\langle K,
H^x\rangle=X$, a contradiction as $\{e\}\not=N_K$ is nilpotent.

\begin{center}
{\bfseries\Large \S~4. Some properties of Carter subgroups}
\end{center}

Here we prove some lemmas that are useful in investigation of Carter subgroups
in finite groups, in particular in almost simple groups.
\vspace{0.5\baselineskip}

\noindent{\bfseries Lemma 4.} {\itshape 
Let $G$ be a finite group satisfying {\em ($*$)}, $H$ be a normal subgroup of
$G$, and $K$ be a Carter subgroup of $G$. Then $KH/H$ is a Carter subgroup
of~${G/H}$.}
\vspace{0.5\baselineskip}

\noindent{\itshape Proof.}
This fact is proven in Lemma 1 under additional assumption
$\vert G\vert\leqslant \vert X\vert$. Now we proved the theorem, thus $\vert
X\vert=\infty$ and this lemma holds true for any finite group~$G$.
\hfill{$\Box$}
\vspace{0.5\baselineskip}

\noindent{\bfseries Lemma 5.} {\itshape
Assume that $G$ is a finite group. Let $K$ be a Carter subgroup of
$G$, with center $Z(K)$. Assume also that  $e\not= z\in Z(K)$ and $C_G(z)$ 
satisfies~{\em($*$)}.
\begin{itemize}
\item[{\em (1)}] Every subgroup $Y$ which contains $K$ and satisfies {\em ($*$)}
is self-normalizing in~$G$.
\item[{\em (2)}] No conjugate of $z$ in $G$, except $z$,
lies in $Z(G)$.
\item[{\em (3)}] If $H$ is a Carter subgroup of $G$, non-conjugate
to $K$, then $z$ is not conjugate to any
element in the center of $H$.
\end{itemize}

In particular the centralizer $C_G(z)$ is self-normalizing
in $G$, and $z$ is not conjugate to any power~${z^k\not= z}$.}
\vspace{0.5\baselineskip}

\noindent{\itshape Proof.}
This lemma is proven in [2, Lemma~3.1] for a minimal counter
example to the problem and therefore its usage for finding Carter subgroups
heavily depends on the classification of
finite simple groups. We state here stronger version of the lemma in order to
avoid such dependence.

(1) Take $x\in N_G(Y)$. Then $K^x$ is a Carter subgroup of $Y$. By the theorem,
Carter subgroups of $Y$ are conjugate. Therefore there exists $y\in Y$ with
$K^x=K^y$. Hence $xy^{-1}\in N_G(K)=K\leq Y$ and $x\in Y$. 

(2) Assume $z^{x^{-1}}\in Z(K)$ for some $x\in G$. Then
$z$ belongs to the center of $\langle G,G^x\rangle\leq C_G(z)$. Since $C_G(z)$
satisfies ($*$), there exists $y\in C_G(z)$ such that $K^x=K^y$. From
$xy^{-1}\in
C_G(z)$, we get $z^{xy^{-1}}=z$ hence $z^x=z^y=z$. We conclude $z^{x^{-1}}=z$.

(3) If our claim is false, substituting $H$
with some conjugate $H^x$ (if necessary), we may assume
$z\in Z(K)\cap Z(H)$, i.~e. $z\in Z(\langle K,H\rangle)\leq C_G(z)$. Again since
$C_G(z)$ satisfies ($*$), there exists $y\in C_G(z)$ such that $H=K^y$. A
contradiction. \hfill{$\Box$}
\vspace{0.5\baselineskip}

Note that for every known finite simple group $G$ (and hence almost simple,
since the group of outer automorphisms is soluble) and for most elements $z\in
 G$ of prime order we have that composition factors of $C_G(z)$ are known
simple groups. Indeed, for sporadic groups this statement can be checked by
using [3]. Composition factors of $C_{A_n}(z)$ are alternating groups.
If $G$ is a finite simple group of Lie type over a field of characteristic $p$
and $(\vert z\vert,p)=1$, then $z$ is semisimple and composition factors of
$C_G(z)$ are finite groups of Lie type. If $\vert z\vert=p$ and $p$ is a good
 prime for $G$, then [4, Theorems~1.2 and~1.4] implies
that all composition factors of $C_G(z)$ are finite groups of Lie type. The
only case, where the structure of centralizers of unipotent elements of order 
$p$ is not
completely known: $p$ is a bad prime for $G$. 
 
Therefore if we are classifying Carter subgroups of almost simple
finite group $A$ by induction we may assume that $C_A(z)$ satisfies ($*$) for
most elements of prime  order $z\in A$. In particular, we can improve table
from [2], using results of present paper and~[5]. In the
table below $A$ is an almost simple group with conjugate Carter subgroups.

{\small $$\begin{array}{|c|c|} \hline \text
{Soc}(A)=G&\text{Conditions for } A\\
\hline
\text{alternating, sporadic};&\\
A_1(r^{t}),\ B_\ell(r^{t}),\
C_\ell(r^{t}),\ t \text{ even if } r=3;&\\
{^2B_2}(2^{2n+1}),\ G_2(r^{t}),\ F_4(r^{t}),\ {^2F_4}(2^{2n+1});&
\text{none}\\
E_7(r^{t}),\ r\not=3;\ E_8(r^{t}),\ r\not=3,5&\\
\hline
D_{2\ell}(r^{t}),\ {^3D_4(r^{3t})},\ {{^2D_{2\ell}}(r^{2t}),}
& A/(A\cap \widehat{G}) \text{ a}\
{2-\text{group}} \\
\ t\ \text{even if }r=3\ {\rm and},\text{ if}\ G=D_4(r^t),& \text{or}\\
|({\rm Field}(G)\cap A):(\widehat
G\cap A)|_{2'}>1 &
|\widehat{G}:(A\cap\widehat G)|\le2\\
\hline B_\ell(3^t),\ C_\ell(3^t),\ D_{2\ell}(3^t),\
{^3D_4(3^{3t})}, \ {{^2D_{2\ell}}(3^{2t})}, &\\
{D_{2\ell+1}(r^{t})},\ {^2D_{2\ell+1}}(r^{2t}),\
{^2G_2(3^{2n+1})}, & A=G\\
E_6(r^{t}),\ {^2E_6}(r^{2t}),\
E_7(3^t),\ E_8(3^t), \ E_8(5^t)&\\
\hline
A_\ell(r^t),\ {^2A_\ell(r^{2t})},\ \ell >1&
G\leq A \leq \widehat{G},\\
\hline
\end{array}$$}

\begin{center}
{\slshape\Large LITERATURE}
\end{center}

\begin{itemize}
\item[{[1]}]
{\slshape F. Dalla Volta, A. Lucchini, M. C. Tamburini},  ``On the Conjugacy
Problem for Carter Subgroups'',  Comm. Algebra, {\bf 26}, No 2, 395--401 (1998).
\item[{[2]}] {\slshape M.C.~Tamburini, E.P.~Vdovin}, ``Carter subgroups of
finite groups'', Journal of Algebra, {\bf 255}, No~1, 148-163 (2002).
\item[{[3]}] {\slshape J.H. Conway, R.T. Curtis, S.P. Norton, R.A. Parker,
R.A. Wilson}, ``Atlas of Finite Groups'', Clarendon Press, Oxford, 1985.
\item[{[4]}] {\slshape G. M. Seitz}, ``Unipotent elements, tilting modules,
and saturation'',  Inv. Math, {\bf 141}, No~3, 467--503 (2000).
\item[{[5]}] {\slshape A.~Previtali, M.~C.~Tamburini, E.~P.~Vdovin},
``The Carter subgroups of some classical groups'',  Bull.LondonMath.Soc.,
{\bfseries 36}, No~1, 145--155 (2004).
\end{itemize}

\end{document}